\documentclass[12pt,a4paper]{article}%
\usepackage{amsfonts}
\usepackage{amsmath}
\usepackage{t1enc}
\usepackage[latin2]{inputenc}
\usepackage[mathscr]{eucal}
\usepackage{amssymb}
\usepackage{amsthm}
\usepackage{graphics}
\usepackage{geometry}
\usepackage{amsxtra}
\usepackage{color}
\usepackage{hyperref}%
\usepackage{graphicx}
\providecommand{\U}[1]{\protect\rule{.1in}{.1in}}
\begin{document}

\title{Distinguish permutations}
\author{István Szalkai\\University of Pannonia, Hungary\\\textsc{szalkai.istvan@mik.uni-pannon.hu}}
\date{2025.12.19.}
\maketitle

\begin{abstract}
We are looking for integer numbers $g_{j}$ and $x_{j}$ ($j=1,...,n$) such that
the sums\quad$T_{\pi}:=\sum_{j=1}^{n}~g_{j}\cdot x_{\pi\left(  j\right)  }$
are different for all permutations $\pi\in S_{n}$ \textbf{and} $\max\left\{
T_{\pi}:\pi\in S_{n}\right\}  $ is as small as possible.

\end{abstract}

\qquad

\subsection*{Introduction}

The source of our mathematical problem is the following (possibly well-known)
magic trick.

\textbf{Trick 0.\quad}\textit{Call three mediums and give them }\textbf{1},
\textbf{2}\textit{\ and }\textbf{3}\textit{\ nuts, respectively. Put three
small objects }\textbf{a,b,c}\textit{\ (e.g. lipstick, pocket knife, pencil)
on the table and }18\textit{\ more nuts. Before leaving the room, give them
the following instruction. Ask them to choose one of the objects. Then, ask
the person, who chose }\textbf{a}\textit{, to take the }same
many\textit{\ further nuts from the table as he/she was given at the
beginning. The person, who chose }\textbf{b}\textit{, would take }twice
many\textit{\ further nuts from the table as he/she was given at the
beginning. Finally, the person, who chose }\textbf{c}\textit{, would take
}four times\textit{\ many further nuts from the table as he/she was given at
the beginning. Ask the mediums to hide their objects and their nuts (except
remained on the table) in their pockets. }

\textit{After returning into the room, looking at the nuts left on the table,
we can guess which medium has hidden }\textbf{a,b,c}\textit{~.}\medskip

The explanation is clear: the number of hidden nuts (equivalently, remained on
the table) is different for each choice of \textbf{a,b,c}. In other words: the
sums\quad$T_{\pi}:=1\cdot\pi\left(  1\right)  +2\cdot\pi\left(  2\right)
+4\cdot\pi\left(  3\right)  =\sum_{j=1}^{3}~2^{j-1}\cdot\pi\left(  j\right)  $
are all different for each permutation $\pi\in S_{3}$~. In this paper we
attempt to generalize this phenomena for more persons with minimal number of
nuts.\bigskip

\subsection*{Notation and basic result}

\qquad Let $n$ be a fixed natural number, $2\leq n$,

$x_{i}\in\mathbb{N}$ natural numbers $=$ the number of nuts given to the
$i^{\prime}th$ person at the beginning,

$g_{j}\in\mathbb{N}$ natural numbers $=$ the multipliers when choosing the
$j^{\prime}th$ object ($i,j=1,...,n$),

$i=\pi\left(  j\right)  $ (%
$<$%
=%
$>$
$j=\pi^{-1}\left(  i\right)  $)
$<$%
=%
$>$
the person $i$ chose the object $j$\ ,

$\pi\in S_{n}$ any permutation.\quad$\Box$ \bigskip

\textbf{1. Our goal} is to find $x_{i}$ and $g_{j}$ such that the sums
\begin{equation}
T_{\pi}:=\sum\limits_{j=1}^{n}g_{j}\cdot x_{\pi\left(  j\right)  }\text{
,}\label{TpiSpi}%
\end{equation}

\qquad\qquad\qquad\qquad are all different for each $\pi\in S_{n}$ \hfill(*)

\noindent and

$\qquad\qquad\qquad\qquad\max\limits_{\pi\in S_{n}}\left\{  T_{\pi}\right\}
$\quad is minimal. \hfill(**)\bigskip

\qquad

\textbf{2. For simplicity }we assume in this paper that
\begin{equation}
g_{1}<g_{2}<...<g_{n}\label{gi<gj}%
\end{equation}
and
\begin{equation}
x_{i}:=i\text{\quad}\left(  i=1,...,n\right)  ~\text{.}\label{xi=i}%
\end{equation}

Further, for any permutation $\pi\in S_{n}$
\begin{equation}
\pi=\binom{1~~\quad~~~2~~\quad~~~.~~.~~.~~\quad~~n}{\pi_{1}~\pi_{2}%
~.~.~.~\pi_{n}}\label{(ipi)}%
\end{equation}
(i.e. $\pi_{i}=\pi\left(  i\right)  $) \textit{we write} simply
\begin{equation}
\pi=\left\langle \pi_{1},\pi_{2},...,\pi_{n}\right\rangle \text{\quad or\quad
}\pi=\left\langle \pi_{n},...,\pi_{2},\pi_{1}\right\rangle ^{\leftarrow}\text{
,}%
\end{equation}
not to confuse with the cyclic permutation\footnote{$^{)}$ The \textit{cyclic
}permutation $\pi=\left(  a_{1},a_{2},...,a_{n}\right)  $ means $\pi\left(
a_{1}\right)  =a_{2}$ , $\pi\left(  a_{2}\right)  =a_{3}$ , ... , $\pi\left(
a_{n}\right)  =a_{1}$ . So, for e.g.\quad$\pi=\left(  1,3,6,2,5,4\right)
=\binom{1~~2~~3~~4~~5~~6}{3~~5~~6~~1~~4~~2}$ ($\pi\left(  1\right)  =3,$
$\pi\left(  3\right)  =6,$ $\pi\left(  6\right)  =2,$ ...), while
\medskip\newline$\rho=\left\langle 1,3,6,2,5,4\right\rangle =\binom
{1~~2~~3~~4~~5~~6}{1~~3~~6~~2~~5~~4}$ .($\rho\left(  1\right)  =1,$
$\rho\left(  2\right)  =3,$ $\rho\left(  3\right)  =6,$ ...). \medskip}$^{)}$
$\left(  \pi_{1},\pi_{2},...,\pi_{n}\right)  $. \bigskip

In what follows, $\overline{g}$ shortens the sequence $\left(  g_{1}%
,...,g_{n}\right)  $.\bigskip

\qquad

\textbf{Statement 3.}\quad Adding the same real number $y$ to each $x_{i}$ the
property (*) remains true: \medskip

$\qquad\qquad T_{\pi}=\sum\limits_{j=1}^{n}g_{j}\cdot\left(  x_{\pi\left(
j\right)  }+y\right)  =\sum\limits_{j=1}^{n}g_{j}\cdot x_{\pi\left(  j\right)
}+y\cdot\sum\limits_{j=1}^{n}g_{j}$

\noindent since\quad$y\cdot\sum\limits_{j=1}^{n}g_{j}$ \ is constant.\qquad
$\Box$\bigskip

\textbf{Theorem 4.}\quad For any $n\in\mathbb{N}$ there is a set $\overline
{g}=\left(  g_{1},...,g_{n}\right)  $ satisfying (*). \medskip

\textbf{Proof:} Let $x_{i}=i-1$ and $g_{j}=n^{j-1}$. Then the sums $T_{\pi}$
are exactly the natural numbers of $n$ digits in the number system of base $n
$. These are clearly different.\quad(Similar is true for $g_{j}=m^{j-1}$ for
any $m\geq n$.)\qquad$\Box$ \bigskip

\textbf{Remark 5.}\quad In the above construction\quad$\max\limits_{\pi\in
S_{n}}\left\{  T_{\pi}\right\}  =n^{n}-1$ , which construction is highly
redundant, since there are (distinct) $n^{n}$ numbers but only $n!$
permutations\footnote{$^{)}$ by the Stirling-formula $n!\sim\left(  \frac
{n}{e}\right)  ^{n}\cdot\sqrt{2\pi n}<n^{n}$.}$^{)}$ and\quad$\lim
\limits_{n\rightarrow\infty}\dfrac{n^{n}}{n!}=\infty$~.\bigskip

\textbf{Examples 6.}\quad In the original Trick 0\quad$n=3$ and $\overline
{g}=\left(  1,2,2^{2}\right)  <\left(  1,3,3^{2}\right)  $ "works" (satisfies
(*)).\quad On the other hand, for $n=4$ the set \ $\overline{g}=\left(
1,3,3^{2},3^{3}\right)  <\left(  1,4,4^{2},4^{3}\right)  $\quad does
\textit{not} work, since e.g. \medskip

$T_{\left\langle 3,1,2,4\right\rangle ^{\leftarrow}}=\left(  3\right)
\cdot\left[  27\right]  +\left(  1\right)  \cdot\left[  9\right]  +\left(
2\right)  \cdot\left[  3\right]  +\left(  4\right)  \cdot\left[  1\right]
=3\cdot27+1\cdot9+2\cdot3+4\cdot1=100$

\noindent and

$T_{\left\langle 2,4,3,1\right\rangle ^{\leftarrow}}=\left(  2\right)
\cdot\left[  27\right]  +\left(  4\right)  \cdot\left[  9\right]  +\left(
3\right)  \cdot\left[  3\right]  +\left(  1\right)  \cdot\left[  1\right]
=2\cdot27+4\cdot9+3\cdot3+1\cdot1=100$ \medskip

For $n=4$ the sets $\overline{g}=\left(  1,3,3^{2},28\right)  $, $\left(
1,3,8,27\right)  $, even $\left(  1,3,8,26\right)  $ are suitable. \bigskip

\subsection*{Antilexicographic order of permutations}

For further investigation we order the permutations in antilexikographic
order: \bigskip

\textbf{Definition 7.}\quad Let $\succ$ denote the \textbf{antilexicographic
order} on $S_{n}$ as: \medskip

$\pi\succ\rho$\quad if there is a $j_{0}\leq n$ such that $\pi\left(
j\right)  =\rho\left(  j\right)  $ for $n\geq j>j_{0}$ (if $j_{0}<n$) and
$\pi\left(  j_{0}\right)  >\rho\left(  j_{0}\right)  $ ,\quad and let\medskip

$\pi\vartriangleright\rho$\quad if $\pi$ is a \textbf{successor} of $\rho$,
i.e. $\pi\succ\rho$ but there is no $\sigma$ such that $\pi\succ\sigma
\succ\rho$~. $\Box$ \medskip

\textbf{Examples 8.}\quad$\left\langle 5,2,4,3,1\right\rangle ^{\leftarrow
}\succ\left\langle 4,5,2,3,1\right\rangle ^{\leftarrow}$ ($j_{0}=n=5$%
),\quad$\left\langle 4,2,3,5,1\right\rangle ^{\leftarrow}\succ\left\langle
4,2,1,5,3\right\rangle ^{\leftarrow}$ ($j_{0}=3<n$) and $\left\langle
5,2,4,3,1\right\rangle ^{\leftarrow}\vartriangleright\left\langle
5,2,4,1,3\right\rangle ^{\leftarrow}$\ ($j_{0}=4<n$).\medskip

\textbf{Remark 9.}\quad Clearly $\pi\left(  j\right)  $ and $\rho\left(
j\right)  $ for $j_{0}>j$ are uninteresting. Especially, in the case $j_{0}=n$
only $\pi\left(  n\right)  >\rho\left(  n\right)  $ is relevant, moreover
$\pi\left(  n\right)  >\rho\left(  n\right)  $ implies $\pi\succ\rho$.

The case $j_{0}=1$ is impossible, since $\pi$ and $\rho$ both are permutations
on $n$ many elements.

Clearly $\succ$ is a complete (linear) ordering on the finite set $S_{n}%
$~;\quad$id=\left\langle n,...,2,1\right\rangle ^{\leftarrow}$ is the greatest
and $\delta=\left\langle 1,2,...,n\right\rangle ^{\leftarrow}$ is the smallest
permutation according to $\succ$~;\quad and for each $\pi\neq id$ there is a
unique $\rho$ such that $\pi\vartriangleright\rho$ ;\quad and for each
$\pi\neq\delta$ there is a unique $\sigma$ such that $\sigma\vartriangleright
\pi$. \bigskip

\textbf{Theorem 10.}\quad$\pi\vartriangleright\rho$\quad\textit{holds if and
only if}\quad i) through vi) all holds: \medskip

i) \quad there is a $j_{0}$ , $2\leq j_{0}\leq n$ such that $\pi\left(
j\right)  =\rho\left(  j\right)  $ for $n\geq j>j_{0}$ \medskip\newline(in
case $j_{0}=n$ no $j$),\quad and\quad$\pi_{j_{0}}>\rho_{j_{0}}$~,\medskip

ii)\quad$\left\{  \pi_{n},...,\pi_{j_{0}+1}\right\}  =\left\{  \sigma
_{n},...,\sigma_{j_{0}+1}\right\}  $\quad(set equality), \hfill($\lozenge
$)\medskip

iii) $H_{j_{0}}:=\left\{  \pi_{j_{0}},...,\pi_{1}\right\}  =\left\{
\sigma_{j_{0}},...,\sigma_{1}\right\}  $\quad(set equality)\medskip

iv)\quad$\pi_{j_{0}}\in H_{j_{0}}\smallsetminus\left\{  1\right\}  $\quad is
arbitrary, \medskip

v)\quad$\rho_{j_{0}}=\max\left\{  \sigma_{j}\in H_{j_{0}}:\sigma_{j}%
<\pi_{j_{0}}\right\}  $ \ \hfill($\triangledown$)\medskip

vi)\quad$\pi_{j_{0}-1}<\pi_{j_{0}-2}<...<\pi_{1}$\quad and\quad$\sigma
_{j_{0}-1}>\sigma_{j_{0}-2}>...>\sigma_{1}$ . $\qquad\Box$ \hfill
\ ($\triangle$)\bigskip

\textbf{Statement 11.}\quad Using $\overline{g}=\left(  1,g,g^{2}%
,...,g^{n-1}\right)  $\quad$\left(  n\leq g\right)  $\ we have
\begin{equation}
\pi\succ\rho\iff T_{\pi}>T_{\rho}\text{ .}\label{<=>}%
\end{equation}
(Use Theorem 4.) \qquad$\Box$\bigskip

\qquad

\textbf{Idea 12. \ }Our next idea to look for sequences $\overline{g}=\left(
g_{1},...,g_{n}\right)  $ satisfying (\ref{<=>}).\medskip

\noindent Since $S_{n}$ is a finite set, (\ref{<=>}) follows from the below
requirement:
\begin{equation}
\pi\vartriangleright\rho\implies T_{\pi}>T_{\rho}\text{ .}\label{=>>}%
\end{equation}

The \textit{case}\ $j_{0}=2$ is tirival: \quad$j_{0}=2$ =%
$>$
$\pi\left(  2\right)  >\rho\left(  2\right)  $ =%
$>$
$\pi\left(  1\right)  =\rho\left(  2\right)  $ and $\rho\left(  1\right)
=\pi\left(  2\right)  $ since $\pi\left(  j\right)  =\rho\left(  j\right)  $
for $n\geq j>j_{0}=2$.\quad So (\ref{=>>}) requires only \medskip

$g_{2}\cdot\pi\left(  2\right)  +g_{1}\cdot\pi\left(  1\right)  >g_{2}%
\cdot\rho\left(  2\right)  +g_{1}\cdot\rho\left(  1\right)  $ , \qquad
i.e.\medskip

$g_{2}\cdot\left(  \pi\left(  2\right)  -\rho\left(  2\right)  \right)
>g_{1}\cdot\left(  \rho\left(  1\right)  -\pi\left(  1\right)  \right)
=g_{1}\cdot\left(  \pi\left(  2\right)  -\rho\left(  2\right)  \right)  $
,\medskip

\noindent which means, that $g_{2}$ and $g_{1}$ can be arbitrary $g_{2}>g_{1}
$~. \bigskip

\subsection*{Properties of $g_{j_{0}}$}

\qquad

\textbf{Theorem 13.}\quad$\pi\vartriangleright\rho$ and (\ref{=>>}) with
($\lozenge$) implies for $j_{0}$ in Theorem 10 \medskip

$\qquad g_{j_{0}}\cdot\pi_{j_{0}}+...+g_{1}\cdot\pi_{1}>g_{j_{0}}\cdot
\sigma_{j_{0}}+...+g_{1}\cdot\sigma_{1}$

\noindent i.e.

\qquad$g_{j_{0}}\cdot\left(  \pi_{j_{0}}-\sigma_{j_{0}}\right)  >g_{j_{0}%
-1}\cdot\left(  \sigma_{j_{0}-1}-\pi_{j_{0}-1}\right)  +...+g_{1}\cdot\left(
\sigma_{1}-\pi_{1}\right)  $ .\quad$\Box$ \hfill($\heartsuit$)\bigskip

The above results offer an induction algorithm for computing the smallest
$g_{j_{0}}$ for $j_{0}=1,2,...n$ :\medskip

\textbf{Algorithm 14. }\medskip

\noindent i)\qquad let $g_{1}$ and $g_{2}$ be arbitrary, $g_{2}>g_{1}$ ,
\medskip

\noindent ii)\quad choose any $j_{0}$ -element set $H_{j_{0}}\subseteq\left\{
1,2,...,n\right\}  $\quad(this is $\binom{n}{j_{0}}$ many choices), \medskip

\noindent iii)\quad let $\pi_{j_{0}}\in H_{j_{0}}\smallsetminus\left\{
1\right\}  $ arbitrary \quad(this is $j_{0}$ or $j_{0}-1$ many choices
)~\footnote{$^{)}$ Total $\binom{n-1}{j_{0}}\cdot j_{0}+\binom{n-1}{j_{0}%
-1}\cdot\left(  j_{0}-1\right)  $ possibilities.}$^{)}$,\medskip

\noindent iv)\quad construct $\left(  \pi_{j_{0}-1},\pi_{j_{0}-2},...,\pi
_{1}\right)  $\ and $\left(  \rho_{j_{0}},\rho_{j_{0}-1},\rho_{j_{0}%
-2},...,\rho_{1}\right)  $\ using ($\triangledown$) and ($\triangle$),\medskip

\noindent v)\quad($\heartsuit$) gives a lower bound for $g_{j_{0}}$ .
\qquad$\Box$ \bigskip

\qquad\bigskip

\subsection*{Numerical results}

\qquad

Algorithm 14. gives the following results for $n=4$ and $n=5$ : \medskip

for $n=4$ : \quad$\ (g1,...,g4)=(1,2,5,15)$, \quad$\ \max\limits_{\pi\in
S_{n}}\left\{  T_{\pi}\right\}  =80<<4^{4}-1=255$ , \medskip

for $n=5$ : \quad\ $(g1,...,g5)=(1,2,6,23,99)$,\quad$\max\limits_{\pi\in
S_{n}}\left\{  T_{\pi}\right\}  =610<<5^{5}-1=3124$ . \medskip

\subsection*{Reference}

\qquad

[Sz] \textbf{Szalkai István:} \textit{An idea} (in Hungarian), Polygon VII.
(1997), 85-88.

\end{document}